\documentclass[11pt]{article}

\title{\bf Reflexion maps and geometry of surfaces in $\RR^4$}
\author{P.J.Giblin \and S.Janeczko \and M.A.S.Ruas  }

\usepackage{amssymb}
\usepackage{amsmath}

\usepackage{epsfig}
\usepackage{color}

\date{}

\newcommand{\RR}{\mathbb R}
\newcommand{\KK}{$\mathcal K$}
\renewcommand{\ll}{\lambda}

\newcommand{\BB}{\mathcal{B}}

\newtheorem{theorem}{Theorem}[section]
\newtheorem{prop}[theorem]{Proposition}

\newtheorem{exam}[theorem]{Example}
\newtheorem{lemma}[theorem]{Lemma}
\newtheorem{defn}[theorem]{Definition}
\newtheorem{rem}[theorem]{Remark}

\begin{document}
\maketitle

\abstract{

In this article we introduce new affinely invariant points---`special parabolic points'---on the parabolic set of a generic surface $M$ in real 4-space, associated with
symmetries in the 2-parameter family of
reflexions of $M$ in points of itself.  The parabolic set itself is detected in this way, and each arc is given a sign, which changes at the
special points, where the family has an additional degree of symmetry.  Other points of $M$ which are detected by the family of reflexions include  inflexion points of
real and imaginary type, and the first of these is also associated with sign changes on the parabolic set.
We show how to compute the special points globally for the case where $M$ is given in Monge form and give some examples illustrating
the birth of special parbolic points in a 1-parameter family of surfaces. The tool we use from
singularity theory is the contact classification of
certain symmetric maps from the plane to the plane and we give the beginning of this classification, including versal unfoldings which we relate to the geometry of $M$.

\noindent
MR 2010 Classification 52A05, 57R45
}

 \section{Introduction}

In a previous article \cite{G-J1} the first two authors studied families of local reflexion maps on surfaces in $\RR^3$ and their
bifurcation sets, in particular showing that certain special parabolic points, not related to the flat geometry of the surface, are detected by the structure of the corresponding
bifurcation set.  These special parabolic, or $A_2^*$ points, arose also in earlier work on centre symmetry sets of surfaces~\cite{G-Z}.  Although the definition of
the reflexion maps is  local the bifurcation sets could be extended over the whole surface, producing curves connecting the special parabolic points.
In this article we extend some of these results to surfaces in $\RR^4$, again studying local reflexions and bifurcation sets of familites of contact maps.
In the present situation we need to study the contact between two surfaces in $\RR^4$ and this is measured by a map (germ) $\RR^2, 0\to \RR^2, 0$.
The appropriate equivalence relation to measure contact is \KK-equivalence (see~\cite{Montaldi}) and therefore the bifurcation set of a family of
contact maps must be constructed according to this equivalence relation, taking into account the inherent $\mathbb Z_2$-symmetry of the contact maps.

We find new `special parabolic points' on a surface in $\RR^4$, which are of two types, `elliptic' and `hyperbolic', and are in some ways analogues of the
special parabolic points encountered in $\RR^3$; the local structure of the bifurcation sets is also similar to the 3-dimensional case.  For a
surface in $\RR^4$ however there are more special kinds of points and the bifurcation set of our family of contact functions displays different structures at these.
We have not so far found a natural interpretation of a {\em global} bifurcation set, connecting special parabolic points and other points through the hyperbolic
and elliptic regions of the surface.

In \S\ref{s:families} we derive the family of reflexion maps and explain our interpretation of the bifurcation set of such a family.
The abstract classification which we need is given in Theorem~\ref{th:classification} and the application to surfaces in $\RR^4$ occupies the remainder of
\S\ref{s:bifsets}. We find the bifurcation set germ at parabolic points, at the two types of special parabolic  points, and at inflexion points of real and imaginary type. In particular we
show that arcs of the parabolic set between these various special points can be given a sign, which changes in a well-defined way at the special points.  Identifying the
local structure of the bifurcation sets requires that we are able to  check versal unfolding conditions and we give the criteria for these to hold in each case.

The above calculations are done with a surface $M$ in Monge form at the origin.
In \S\ref{s:examples} we show how to compute the special parabolic points on a whole surface given in Monge form.  The special parabolic points
are found as the intersection of the parabolic set with another curve in $M$ and we find an explicit formula for this curve, given in Appendix A but applied
to several examples in \S\ref{s:examples}. An example, adapted from~\cite{B-T}, shows the birth of special parabolic points on a loop of the parabolic set
created in a generic 1-parameter family of surfaces---an elliptic island in a hyperbolic sea.  Immediately after the moment that the island appears it has
no special parabolic points but two of these, of the same type, can be born as the island grows larger.  Between the two the sign of the parabolic set changes.

Finally in \S\ref{s:qu} we give some concluding remarks and open problems.

\section{Families of contact maps}\label{s:families}
Consider a surface $M$ in $\RR^4$, with coordinates $(a,b,c,d)$, parametrized by $\gamma(x,y)= (f(x,y), g(x,y),x,y)$, where we shall assume that
the 1-jets of $f$ and $g$ at $(x,y)=(0,0)$ are zero. Let $(p,q)$ be the parameters of
a fixed point on the surface. Reflecting a point $\gamma(p+x,q+y)$ of $M$ in the point $\gamma(p,q)$ gives $2\gamma(p,q)-\gamma(p+x, q+y)$, so that
reflecting $M$ in $\gamma(p,q)$  gives the surface $M^*$ through $\gamma(p,q)$ with parametrization $\RR^2\to\RR^4$:
\[ (x,y)\mapsto (2f(p,q)-f(p+x,q+y), \  2g(p,q) - g(p+x,q+y), \ p-x, \ q-y).\]
Thus $x=y=0$  returns the point $\gamma(p,q)$.
Composing this parametrization with the map $\RR^4\to\RR^2$ defined by $(a,b,c,d)\mapsto (f(c,d)-a, \ g(c,d)-b)$, for which the inverse image of $(0,0)$
is  equal to $M$, gives the following  map (germ) $F_{(p,q)}: \RR^2, (0,0) \to \RR^2, (0,0)$, whose \KK-class measures the
contact  between $M$ and $M^*$ at $\gamma(p,q)$ (see~\cite{Montaldi}).
\begin{eqnarray}
F_{(p,q)}(x,y)&=& (f(p+x,q+y)+f(p-x,q-y)-2f(p,q), \nonumber\\
&& g(p+x,q+y)+g(p-x,q-y)-2g(p,q)).
\label{eq:contact-family}
\end{eqnarray}
When we include the parameters $p,q$ we write $F(x,y,p,q)$.  Note that $F(x,y,p,q)\equiv F(-x,-y,p,q)$:
for each $(p,q)$ the map $F_{(p,q)}$ is symmetric with respect to the reflexion $(x, y) \to (-x,-y)$.

Thus $F$ is a family of symmetric mappings $\RR^2\to \RR^2$, with variables $x,y$ parametrized by $p,q$. We  investigate the bifurcation set of this family, the fundamental
definition of which is
\begin{eqnarray*}
\mathcal B_F&=&\{ (p,q) : \mbox{ there exist } (x,y) \mbox{ such that } F_{(p,q)}  \mbox{ has an  unstable} \\ && \mbox{ singularity at } x,y 
 \mbox{ with respect to } \mathcal K 
 \mbox{ equivalence} \\ && \mbox{ of maps symmetric in the above sense}\}.
\end{eqnarray*}

In \cite{G-J1} the corresponding bifurcation set of a family $F$ of real-valued {\em functions}
was analysed by studying the critical set of $F$.   Here we need to work directly with \KK-equivalence of maps, where the
critical set does not play so significant a role, and we adopt a different approach.

At $(p,q)=(0,0)$ the contact map is
\begin{equation}
F_{(0,0)}(x,y) = (f(x,y) + f(-x,-y), \ g(x,y) + g(-x,-y)),
\label{eq:contact-at-0}
\end{equation}
which is twice the even part of $(f,g)$, but we shall sometimes ignore the factor 2.
Thus the conditions on $M$ needed for the classification of $F_{(0,0)}$ involve only the {\em even
degree} terms of $f,g$; however the conditions for the family $F$ with parameters $p,q$ to give a \KK-versal unfolding will involve also the odd degree terms.

We work within the set of maps $h: \RR^2\to\RR^2$ which are symmetric by reflexion in the origin: $h(x,y)=h(-x,-y)$. To do this we use the basis $u=x^2, v=xy, w=y^2$
for all functions of two variables which are symmetric with respect to this symmetry and study map germs
$H: \RR^3 \to \RR^2$ with coordinates $(u,v,w)$ in $\RR^3$, up to \KK-equivalence preserving the
homogeneous variety (cone) $V: v^2=uw$. (In fact for us this is a half-cone since $u=x^2$ and $w=y^2$ are non-negative, but for classification
purposes we may assume that the whole cone is preserved.) We write
$_V\mathcal{K}$-equivalence for this equivalence of  germs $H: \RR^3, (0,0,0) \to \RR^2,(0,0)$. We shall work with $_V\mathcal{K}$-versal unfoldings and construct
bifurcation diagrams for these in a sense we now explain.

For a given germ $H$, the $_V\mathcal{K}$ equivalence will preserve the intersection $H^{-1}(0) \cap V$ up to local diffeomorphism of $\RR^3$, and indeed
will preserve the multiplicity of intersection of the curve $H^{-1}(0)$ with the cone $V$. As the map $H$ varies in a family the multiplicity will change
and furthermore intersection points of multiplicity $>1$ may move away from the origin; these points nevertheless form part of the `contact data' of $H^{-1}(0)$ and $V$ since
they represent unstable mappings.  Except in one case, described  below, all the contact data are concentrated at the origin.

\begin{defn} \label{def:bifset}
The strata of our bifurcation set are those points in the versal unfolding space for which the contact data consisting of the multiplicity of contact between
$H^{-1}(0,0)$ and $V$ in an arbitrarily small neighbourhood of the origin in $\RR^3$ are constant.
\end{defn}

\begin{figure}[!ht]
\begin{center}
\includegraphics[width=3.7in]{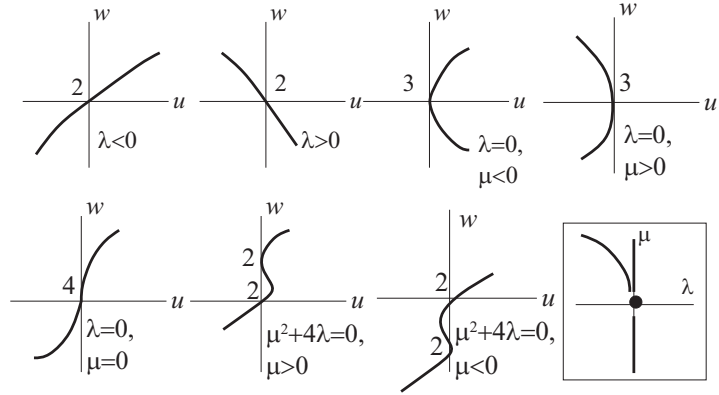}
\end{center}
\caption{\small  The unstable intersections between the  curve $v=0, \ u=w^3-\ll w -\mu w^2$ and the cone $V: v^2=uw$ for various values of $\lambda, \mu$.  These give 0- and 1-dimensional strata
of the bifurcation set of the family $H(u,v,w)=(v, u-w^3+\mu w + \ll w^2)$, shown in the boxed diagram at bottom right.  Intersections corresponding to real values
of $(x,y)$ are in the quadrant $u\ge 0, w\ge 0$ of the $u,w$-plane. }
\label{fig:special-parab-uw}
\end{figure}

The idea
is best illustrated by an example, which will arise in \S\ref{ss:specialparabolic} below.
Consider the family of maps $H_{\ll,\mu}(u,v,w) = (v, u - w^3 + \ll w + \mu w^2)$. For any $(\ll, \mu)$, $H_{\ll,\mu}^{-1}(0)$ lies in the plane $v=0$ with coordinates $(u,w)$,
and $V: v^2 = uw$ intersects this plane in the two lines $u=0, w=0$ (for real solutions for $x,y$ we require indeed $u\ge 0$ and $w\ge 0$).
We therefore examine how the curve $u-w^3+\ll w + \mu w^2=0$ in the $(u,w)$ plane
meets the two coordinate axes.  Intersection with the axis $w=0$ gives only the origin.
Intersection with the axis $u=0$ requires $w(-w^2+\mu w + \ll)=0$ which gives tangency at the origin when $\ll=0$, so that in the $(\lambda, \mu)$ plane
the axis $\ll=0$, apart from the origin,
is one stratum of the bifurcation set. The total contact between $H_{\ll,\mu}^{-1}(0,0)$ and $V$ at the origin is 3.
The origin $\ll=\mu=0$ is a separate stratum since the contact there between $H_{\ll,\mu}^{-1}(0,0)$ and $V$ is 4. There is also a double root
of $-w^2+\mu w + \ll=0$ at $w=\frac{1}{2}\mu$
when $\mu^2 + 4\ll = 0$, resulting in ordinary tangency between $H_{\ll,\mu}^{-1}(0,0)$ and $V$ at $(u,w)=(0,\frac{1}{2}\mu)$.
This gives a stratum $\mu^2+4\ll=0$  of the bifurcation set, with $\mu \ge 0$ since $w=y^2\ge 0$, which intersects every neighbourhood of $(0,0)$ in the plane of the
unfolding parameters $(\ll, \mu)$. The various possibilities are sketched in Figure~\ref{fig:special-parab-uw}  where the intersection number between $H_{\ll, \mu}=0$ and $V$ is
indicated against each intersection point. For real solutions $(x,y)$ we require these intersection points to be in the  quadrant $u\ge 0, w \ge 0$ of the $(u,w)$ plane.  The resulting bifurcation set
is also drawn in Figure~\ref{fig:special-parab-uw}, with four strata of positive codimension in the $(\ll,\mu)$ plane.

\section{Classification of the contact maps up to $_V\mathcal{K}$-equivalence} \label{s:bifsets}

 We consider map germs $\RR^3\to \RR^2$, with coordinates $u,v,w$ in the source ($u=x^2,v=xy,w=y^2$ as above),
under contact equivalence which preserves the homogeneous variety $V: uw-v^2=0$.  Vector fields generating those tangent
to this variety are given by the Euler vector field and the three hamiltonian vector fields:
\begin{equation} u\textstyle{\frac{\partial}{\partial u}}\displaystyle +  v\textstyle{\frac{\partial}{\partial v}}\displaystyle +  w\textstyle{\frac{\partial}{\partial w}}\displaystyle, \
 2v\textstyle{\frac{\partial}{\partial u}}\displaystyle +  w\textstyle{\frac{\partial}{\partial v}}\displaystyle, \
 u\textstyle{\frac{\partial}{\partial v}}\displaystyle +  2v\textstyle{\frac{\partial}{\partial w}}\displaystyle, \
 u\textstyle{\frac{\partial}{\partial u}}\displaystyle -  w\textstyle{\frac{\partial}{\partial w}}\displaystyle.
\label{eq:vfs}
\end{equation}
The tangent space to the $_V\mathcal K$ orbit at $H(u,v,w)$ is
$dH(\theta_V) + H^*(m_2)\mathcal{E}_3^2,$
where $dH$ is the jacobian matrix of $H$ and $\theta_V$ is the $\mathcal E_3$ module generated by the above vector fields.

\medskip

The classification which we need is summarized in Theorem~\ref{th:classification}, which is proved by the method of complete transversals \cite{BKdP}
and the finite determinacy theorem for $_V\mathcal K$ equivalence \cite{D}. Comments on this
classification and  application to our geometrical situation are in the remainder of this section. (We remark here that a different but related classification of
maps involving only odd degree terms is obtained in \cite{MRD}.)
\begin{theorem}
The abstract classification of map germs $H:\RR^3\to\RR^2$ up to \KK-equivalence preserving the half-cone $V:v^2-uw=0, \ u\ge 0, \ w\ge0$ starts with the  classes given in Table~\ref{table}.
The classes of symmetric germs $h:\RR^2\to\RR^2$, where $h(x,y)=h(-x,-y)$, up to \KK-equivalence preserving the symmetry are obtained by replacing
$u,v,w$ by $x^2, xy, y^2$ respectively. \hfill$\Box$
\begin{table}[!ht]
\begin{center}
\hspace*{-2cm}
\begin{tabular}{|c|c|c|c|c|} \hline
type &normal form & $_V\mathcal{K}$ codimension & versal unfolding  & geometry \\ \hline\hline
(H) & $(w,u)$  & 0 &  --- & hyperbolic point \\ \hline
(E) & $(u-w, v)$ & 0 & --- &  elliptic point \\ \hline
(P) & $(v, u\pm w^2)$ & 1 & $(0, \ll w)$ &ordinary parabolic point \\ \hline
(SP) & $(v, u\pm w^3)$ & 2 & $(0,\ll w + \mu w^2)$ & special parabolic point \\ \hline
(IR) & $(v, u^2 + 2buw \pm w^2)$& 3 & $(0, buw+ \lambda u + \mu w)$  & inflexion of real type \\
& $b^2\ne 1$ for $+$ &  && \\ \hline
(II) & $(u+w, k u^2 + uv)$ & 3 & $(0, ku^2+\lambda u + \mu v)$& inflexion of imaginary type \\
& or $(u+w, uv + kv^2)$ &3 &$(0, kv^2 + \lambda u + \mu v)$ & \\ \hline
\end{tabular}
\end{center}
\vspace*{-0.2in}
\caption{The lowest codimension singularities in the $_V\mathcal K$ classification of map germs $\RR^3,(0,0,0)\to\RR^2,(0,0)$.}
\label{table}
\end{table}
\label{th:classification}
\end{theorem}
We shall see that the moduli $b$ and $k$ in the normal forms above do not affect the geometry of the situation. Note that the two
forms $(v,u\pm w^2)$ are not equivalent since $u\ge 0$ so we cannot replace $u$ by $-u$.  The same applies to the two forms
$(v, u\pm w^3)$. Note that the germs (P) and (SP) are the first two in a sequence $(v, u\pm w^k), \ k\ge 2$, distinguished by the contact
between the zero-set of the germ and the cone $V: v^2-uw=0$.

\medskip

 The contact maps are invariant
under affine transformations of the space $\RR^4$ in which our surface $M$ lies, so that we may first put $M$ in a standard
form at the origin in $(a,b,c,d)$-space. We can assume the tangent plane at the origin is the $(c,d)$-plane and the
quadratic terms $f_2, \ g_2$ of $f,\ g$ are reduced
by the action of $GL(2,\RR)\times GL(2,\RR)$ on pairs of binary quadratic forms to a standard form. Finally
a linear transformation of $\RR^4$ reparametrizes $M$ as $(x,y)\mapsto (f, g, x, y)$ where now
$f$ and $g$ have their quadratic parts in standard form. See for example
\cite[pp. 182--183]{B-T} for the standard forms of 2-jets of surfaces in $\RR^4$.

There is a convenient way to recognize the types (P) and (SP) of the contact map $(u,v,w)\mapsto (C_1(u,v,w), C_2(u,v,w))$, which will be useful below.
\begin{lemma}
In each case the zero set $C_1=C_2=0$ in $\RR^3$ is a smooth curve at the origin
and \\
 (P): has exactly 2-point contact (ordinary tangency) with the cone $V: v^2=uw$ at the origin,\\
(SP): has exactly 3-point contact with the cone $V$ at the origin. \hfill$\Box$
\label{lemma:PSP}
\end{lemma}

\subsection{First stable case: hyperbolic point}\label{ss:hyperbolic}  A standard form
for the  2-jet of the surface at a hyperbolic point is $(y^2, \ x^2, \ x, \ y)$, or in a less reduced form
 $(f_{11}xy + f_{02}y^2, \ g_{20}x^2, \ x, \ y)$ where $f_{02} \ne 0, g_{20}\ne 0$.
The contact map at the origin of $\RR^4$, ignoring the factor 2 in (\ref{eq:contact-at-0}), has 1-jet $F_1 = (f_{11}v+f_{02}w, \ g_{20}u)$
(or just $(w,u)$ in the reduced form).
This is $_V\mathcal{K}$-stable and is the case where
the kernel of the linear map $F_1: \RR^3\to\RR^2$
intersects the cone $V\subset \RR^3$ only in the origin. The bifurcation set germ is empty.

\subsection{Second stable case: elliptic point}\label{ss:elliptic}  A standard form
for the  2-jet of the surface is $(x^2-y^2, \ xy, \ x, \ y)$, or in a less reduced form
 $(f_{20}x^2 + f_{02}y^2, \ g_{11}xy, x, y), \ f_{20}f_{02}<0, g_{11}\ne 0$
as in~\cite{B-T}.
This corresponds to 1-jet  $F_1= (f_{20}u+f_{02}w, \ g_{11}v)$ (or $(u-w, \ v)$ in reduced form) .
This is $_V\mathcal{K}$-stable and it is the case where
the kernel of the linear map $F_1 : \RR^3\to\RR^2$ intersects the cone $V$ in two distinct generators.
 The bifurcation set germ is empty.

\subsection{Codimension 1 case: ordinary parabolic point}\label{ss:parabolic}

A standard form of the 2-jet of $M$ at a parabolic point, up to affine transformations of $\RR^4$, is
$(f_{11}xy, \ g_{20}x^2, \ $ $x, \ y)$, where $ f_{11}\ne 0, g_{20} \ne 0$. The corresponding 1-jet in $(u,v,w)$ coordinates is
$(v, u)$ from the abstract classification, with gives 2-jet $(v, u\pm w^2)$ which is 2-$_V\mathcal K$-determined.
The two cases, with signs $\pm$, are not equivalent. Note that with 1-jet $(v,u)$ the kernel of the linear
map from $\RR^3$ to $\RR^2,\ (u,v,w)\mapsto (v,u)$, is  along the $w$-axis, which is a generator of the cone $V$.

For the contact map $F_{(0,0)}$ we obtain $(f_{11}xy, g_{20}x^2\pm g_{04}y^4)$,
provided the coefficient $g_{04}$ of $y^4$ is nonzero, with two cases according as $g_{20}g_{04}>0$ or $<0$.
(It can be checked that in reducing to this form the coefficients of $(0,x^2)$ and $(0,y^4)$ are not changed, in particular the
final values are not influenced by the coefficients in the polynomial $f$, provided of course that $f_{11}\ne 0$.)
The coefficient of $y^2$ in the expansion of the second component of $F(x,y,p,q)$ is $2g_{12}p + 6g_{03}q$; thus provided
$g_{12}\ne 0$ or $g_{03}\ne 0$ the family~(\ref{eq:contact-family}) with parameters $p, q$ gives
a versal unfolding (note that these are odd degree terms of $g(x,y)$).
We call such points, where the expansion of $M$  at the origin has the 2-jet $(f_{11}xy, \ g_{20}x^2, \ x, \ y)$ and
\begin{equation}
f_{11} \ne 0, \ g_{20} \ne 0, \ g_{04} \ne 0, \ g_{12} \mbox{ or } g_{03} \ne 0,
\label{eq:ordparab}
\end{equation}
{\em ordinary parabolic points} of $M$. The last condition is equivalent to the smoothness of the
parabolic set of $M$ at the origin (see below) but the condition $g_{04}\ne 0$ does not arise from the
flat geometry of $M$ and is analogous to the condition found in~\cite{G-J1} for an `ordinary' ($A_2$) point
of the parabolic set of $M\subset \RR^3$.

A standard result is that the global equation of the
parabolic set of a surface $M$ in the form $(f(x,y), g(x,y), x, y)$ is
\begin{equation}
 (f_{xx}g_{yy}-f_{yy}g_{xx})^2=4(f_{xy}g_{yy}-f_{yy}g_{xy})(f_{xx}g_{xy}-f_{xy}g_{xx}).
\label{eq:parab-locus}
\end{equation}
This can be proved by considering the 3-parameter family of height functions at any point of $M$, say
$H(x,y,\lambda, \mu, \nu)=\lambda f(x,y)+g(x,y)+\mu x+\nu y$ or
$H(x,y,\lambda, \mu, \nu)= f(x,y)+\lambda g(x,y)+\mu x+\nu y$
and writing down the condition that there is a unique  direction $(\lambda, 1, \mu, \nu)$ or
$(1, \lambda, \mu, \nu)$ with the height function having
a non-Morse singularity, that is $H_x=H_y=H_{xx}H_{yy}-H_{xy}^2=0$. (All normal vectors to
$M$ have one of these two forms.) We note below in \S\ref{ss:formulas} that
the formula also follows from our analysis of contact functions.

\medskip

In the present case the lowest terms in the equation of the parabolic set at the origin are, from (\ref{eq:parab-locus}),
 $16f_{11}^2g_{20}(g_{12}x + 3g_{03}y)$,
so that the parabolic set is smooth at the origin if and only if $g_{12}$ or $g_{03}$ is nonzero: the last condition of~(\ref{eq:ordparab}).
We can unambiguously label smooth segments of the parabolic set with the sign $+$ or $-$ according as, with 2-jet
of $(f,g)$ equal to $(f_{11}xy, \ g_{20}x^2)$, both coefficients being nonzero, the product $g_{20}g_{04}$ of the coefficients
of $(0,x^2)$ and $(0,y^4)$ is $>0$ or $<0$. We shall see below when the sign of the parabolic set changes.

For the bifurcation set, we consider the map $(u,v,w)\mapsto (v, \ u \pm w^2 + \lambda w)$ and the multiplicity of the zero set of this in
an arbitrarily small neighbourhood of the origin.
Since $v=0$ the intersection lies in the $(u,w)$ plane, at points of the $u$- and $w$-axes. The curve $u=\mp w^2-\ll w$ is tangent
to the $w$ axis if and only if $\ll=0$ and then the multiple value of $w$ is 0 so the tangency is at the origin.
In the geometrical case of a surface, as above, the condition $\lambda=0$ is replaced by $2g_{12}p + 6g_{03}q=0$,  which is the tangent line
to the parabolic set at the origin.  Thus the germ of the bifurcation diagram in the $(p,q)$ parametrization plane of the surface consists of the
tangent line to the parabolic set:
\begin{prop} At an  point of the parabolic set satisfying (\ref{eq:ordparab}) the bifurcation set $\BB$ is locally exactly the parabolic set. We can give a sign
to each such point of the parabolic set by the sign of $g_{20}g_{04}$ when the $2$-jet of $(f,g)$ is reduced to $(f_{11}xy, \ g_{20}x^2)$.
\label{prop:parab}
\end{prop}
Points off the parabolic set have stable contact maps,
 in fact they are elliptic or hyperbolic points as in \S\S\ref{ss:hyperbolic} and~\ref{ss:elliptic}.

\subsection{Formulas for loci of types (P) and (SP) in Table~\ref{th:classification}}\label{ss:formulas}

We can use the criterion in Lemma~\ref{lemma:PSP} to obtain the equation (\ref{eq:parab-locus}) for the parabolic set on a general surface in Monge form,
and then find an additional equation which holds at special parabolic points.  We shall use these in \S\ref{s:examples} to analyse some examples of special
parabolic points.

For the contact map (\ref{eq:contact-family}) at the point of $M$ with parameters $p,q$
write $f_{11}$ for  $f_{xx}(p,q)$, $f_{12}$ for $f_{xy}(p,q)$, $f_{1222}$ for $f_{xyyy}(p,q)$ and so on.  Then the 2-jet of the
first component of the contact map $F=F_{(p,q)}$ in terms of $u,v,w$ is (taking into account the factor 2 which automatically arises)
\[ C_1(u,v,w)=\left(f_{11}u + 2f_{12}v + f_{22}w\right) +\]
\[ \textstyle{\frac{1}{12}}\displaystyle\left(f_{1111}u^2 + 4f_{1112}uv +6f_{1122}uw + 4f_{1222}vw + f_{2222}w^2\right),\]
with a similar formula for the second component.

We can now solve the equations $C_1=C_2=0$ for say $u$ and $v$ in terms of $w$ up to order 2, and substitute in the equation $v^2=uw$ of the cone $V$ to
obtain the order of contact of the zero set of $C$ with $V$. The condition for
the order of contact to be at least 2, that is the condition for the coefficient of $w^2$ after substitution to be zero,
then works out at exactly (\ref{eq:parab-locus}) where $f_{xx}$ appears as $f_{11}$ and so on.

The additional condition for the contact to be of order at least 3, that is for the coefficient of $w^3$ also to be zero,
 is naturally more complicated and requires solving for $u$ and $v$ as above to a higher order. But it is possible
to use this condition in explicit examples and it is stated in appendix~\ref{s:app}.  This formula is used in examples in \S\ref{s:examples}.

\subsection{Codimension 2 case: special parabolic point}\label{ss:specialparabolic}
\begin{figure}[!ht]
\begin{center}
\includegraphics[width=4.5in]{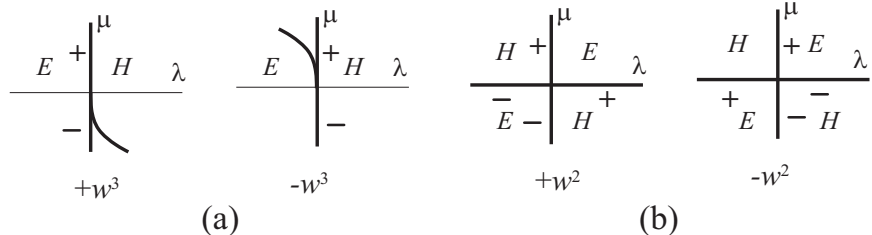}
\end{center}
\vspace*{-0.3in}
\caption{\small (a) The bifurcation set of the unfolding $(v, \ u\pm w^3 +\lambda w + \mu w^2)=(xy, \ x^2 \pm y^6 + \lambda y^2 + \mu y^4)$,
as in \S\ref{ss:specialparabolic} (special parabolic points), with $+w^3$ on the left and $-w^3$ on the right. The bifurcation set in each case consists of a germ of the $\mu$-axis and a half
parabola.  In the geometrical situation the $\mu$-axis corresponds to the parabolic set of $M$ and the sign, $+$ or $-$,
against this axis is the sign attached to that segment of the parabolic set as in \S\ref{ss:parabolic}. Further $E$ and $H$
refer to the parts of the $(\lambda, \mu)$ plane which correspond with elliptic and hyperbolic points of $M$, respectively,
using the normal forms of \S\S\ref{ss:hyperbolic},\ref{ss:elliptic}. The left-hand figure of (a) corresponds with $4g_{20}g_{06}-g_{13}^2>0$
and the right-hand figure with $4g_{20}g_{06}-g_{13}^2<0$. \newline
 (b) Similarly for the bifurcation set of $(v, u^2 \pm w^2 +\lambda u + \mu w)$
as in \S\ref{ss:realinflexion} (inflexions of real type), corresponding to  $4g_{40}g_{04}-g_{22}^2 > 0$ on the left and $<0$ on the right in the geometrical situation. }
\label{fig:bifset-specialparabolic-realinflexion}
\end{figure}

This degeneracy occurs for the abstract map $\RR^3 \to \RR^2$ when the coefficient of $w^2$ in~\S\ref{ss:parabolic} equals zero but there is a nonzero
coefficient of $w^3$. The kernel of the 1-jet map $\RR^3\to\RR^2, \ (u,v,w) \mapsto (v,u)$ is still 1-dimensional
and along a generator of the cone $V$. The bifurcation set of the abstract germ in the case $(v, u-w^3)$ was analysed in \S\ref{s:families} and is illustrated in Figure~\ref{fig:special-parab-uw}.
The other case, $(v,u+w^3)$, is similar and the full picture of the bifurcation set is in Figure~\ref{fig:bifset-specialparabolic-realinflexion}.

In our geometrical situation, on the surface $M$ the above degeneracy
corresponds to a parabolic point with the 2-jet of $(f,g)$ being $(f_{11}xy, \ g_{20}x^2) $
and $g_{04} = 0$. The additional condition which ensures that the contact singularity is no more degenerate
 is $g_{13}^2 - 4g_{20}g_{06}\ne 0$, that is the even degree terms $g_{20}x^2+g_{13}xy^3+g_{06}y^6$ do not form a perfect square.
(This condition remains unchanged when the higher terms of $f$ are eliminated, in particular the condition to avoid
further degeneracy does not involve the higher degree terms of $f$.)
 We call these {\em special parabolic points}\footnote{In the case of a surface in $\RR^3$ they had an alternative name,
``$A_2^*$ points'', referring to the fact that the contact between the surface and its tangent plane at any parabolic point is
a function of type $A_2$, but
this notation is not appropriate here.}.
The further condition that in the family of contact maps the parameters $p,q$ give a versal unfolding is $5g_{12}g_{05}-3g_{03}g_{14}\ne 0$.

The two cases, distinguished by the sign of $g_{13}^2-4g_{20}g_{06}$ in the geometrical
situation, differ as to the region of $M$, elliptic or hyperbolic, in which the `half parabola' branch of $\BB$ lies.
Figure~\ref{fig:bifset-specialparabolic-realinflexion}(a) shows the two cases. Furthermore, at points along the parabolic set, the local expansion of the surface has
$g_{04} \ne 0$ and $g_{04}$ changes sign at special parabolic points. Thus if we label points of the parabolic set by $+$ or $-$
then the sign changes at special parabolic points. See Figure~\ref{fig:bifset-specialparabolic-realinflexion}(a).

Summing up the conclusions of this section:
\begin{prop}
A parabolic point of $M$, with the $2$-jet of $(f,g)$ in the form $(f_{11}xy, \ g_{20}x^2) $
is called a {\em special parabolic point}  if the coefficient $g_{04}$ of $y^4$ in $g$
is zero and $g_{13}^2 - 4g_{20}g_{06}\ne 0$. The sign attached to ordinary parabolic points close
to this one, as in Proposition~\ref{prop:parab}, changes at a special parabolic point. Provided $5g_{12}g_{05}-3g_{03}g_{14}\ne 0$
the $p,q$ parameters versally unfold the contact singularity in the family $F$ as in~(\ref{eq:contact-family}) and the bifurcation set is the union of
the parabolic set and a ``half-parabola'' lying in the  hyperbolic or elliptic region according to the sign of
 $ 4g_{20}g_{06}-g_{13}^2 $, as in Figure~\ref{fig:bifset-specialparabolic-realinflexion}(a).
\end{prop}
We do not know whether there is any significance
attached to the sign of $5g_{12}g_{05}-3g_{03}g_{14}$.

\subsection{First codimension 3 case: inflexions of real type}\label{ss:realinflexion}
The 2-jet of $(f,g)$ at inflexion points of real type (also called real inflexions or umbilic points) on $M$ has the form $(f_{11}xy, \ 0)$,
where $f_{11}\ne 0$.

The abstract map $\RR^3\to\RR^2, \ (u,v,w)\mapsto (v,0)$ has a 2-dimensional kernel which intersects the cone $V$ along two generators. The
abstract normal form is $(v, u^2+2buw \pm w^2)$ where the second component should not be a perfect square, that is $b^2\ne \pm1$
(for the $-$ case this is no restriction).
An abstract $_V\mathcal{K}$-versal
unfolding is given by $(v, u^2 +2buw \pm w^2 +\lambda u + \mu w)$, that is $b$ is a smooth modulus in this case. The bifurcation set $\BB$ is found by considering
the contact of the curve  $u^2 +2buw \pm w^2 +\lambda u + \mu w=0$ with the $u$ and $w$ axes in the $(u,w)$ plane. The condition for tangency comes
to $\mu=0$ or $\lambda=0$, irrespective of the sign in the normal form and the value of $b$. Thus $\BB$ consists of the complete
$\lambda$ and $\mu$ axes (not half-axes), and does not depend on the modulus $b$.
Note that although $uw=v^2$ on the cone $V$ our map germs are defined on $\RR^3$ and not just on the cone, so we cannot use left-equivalence
to remove the modulus term $2buw$.

\begin{rem}
{\rm
 We do not know if $b$ has any geometrical significance. However, taking the two components of the map $(v, u^2+2buw \pm w^2)$, the intersection of the cone $V$ with the plane $v=0$
gives two lines in the plane, $u=0$ and $w=0$, and the second component gives two more lines which are real when $b^2 > \pm 1$ (no restriction for the $-$ sign). The
cross-ratio of these four lines will be responsible for the existence of a smooth modulus.
}\label{rem:mod1}
\end{rem}

The contact singularity for $\lambda=0, \mu\ne 0$
or $\mu=0, \lambda\ne 0$ is equivalent to that for a parabolic point as in~\S\ref{ss:parabolic}.  Thus the two crossing
branches of $\BB$ represent, in our geometrical situation, the parabolic set on $M$.  Indeed at a generic inflexion of real
type the parabolic set does have a transverse self-crossing. Furthermore, as $\lambda$ passes through zero
the normal form for the contact singularity at a parabolic point changes from the $+$ case to the $-$ case or vice versa;
similarly when $\mu$ passes through zero. So the sign attached to the parabolic set changes along each branch
of $\BB$ at an inflexion point of real type.

In the geometrical situation, on the surface $M$ the bifurcation set divides the surface locally into
four regions, two opposite regions being hyperbolic and two elliptic. The configuration corresponding to the two
normal forms is shown in Figure~\ref{fig:bifset-specialparabolic-realinflexion}(b).
The condition to avoid further degeneracy is $g_{22}^2-4g_{40}g_{04} \ne 0$ and the condition for
$p,q$ in the family of contact maps to versally unfold the singularity is $9g_{30}g_{03}-g_{12}g_{21}\ne 0$. Perhaps surprisingly, this
latter condition is the same as that for an inflexion point of real type to be $\mathcal R^+$ versally unfolded by the family of height functions.
(See\footnote{Translating notation from this to our notation we have $a_{20}=f_{20}=0, a_{02}=f_{02}=0,
a_{21}=f_{21}, b_{30}=g_{30}, b_{31}=g_{21}, b_{32}=g_{12}, b_{33}=g_{03}$.  The condition in \cite{IRRT} for a versally unfolded $D_4$ then
reduces to our $\frac{1}{2}f_{11}(-9g_{30}g_{03}+g_{21}g_{12})\ne 0.$  Of course we do not have a $D_4$ singularity, that is the nondegeneracy
of the degree 3 terms of $g$ does not apply.  Instead we have a nondegeneracy condition on the degree 4 terms of $g$.} \cite[Prop.7.9, p.224]{IRRT}.)
As above, the bifurcation set consists of the two intersecting branches of the parabolic set, and passing through the
crossing point on either branch the ``sign'' of the parabolic set, as in \S\ref{ss:parabolic}, changes. See Figure~\ref{fig:bifset-specialparabolic-realinflexion}.

\begin{prop}
At a generic inflexion point of real type on $M$ the $_V\mathcal{K}$ bifurcation set of the family of contact maps consists of the two branches of
the parabolic set through the inflexion point.  The sign as in \S\ref{ss:parabolic}  changes along each branch. See
Figure~\ref{fig:bifset-specialparabolic-realinflexion}(b).
With $2$-jet of $(f,g)$ equal to  $(f_{11}xy, \ 0)$, where $f_{11}\ne 0$, the conditions are $g_{22}^2-4g_{40}g_{04} \ne 0$ and
$9g_{30}g_{03}-g_{12}g_{21}\ne 0$.
\end{prop}

\subsection{Second codimension 3 case: inflexion point of imaginary type}\label{imaginaryinflexion}
The 2-jet of $(f,g)$ at inflexion points of imaginary type on $M$ (also called imaginary inflexions or umbilic points) has the form $(f_{20}x^2 +f_{02}y^2, \ 0)$,
where $f_{20}f_{02}>0$.

The abstract map $\RR^3\to\RR^2$ has kernel of the linear part $(u+w,\ 0)$, a plane meeting the cone $V$ only in the origin, and
reduces to the abstract normal form $H(u,v,w)=(u+w, \   au^2 + 2buv + cv^2)$,
subject to the conditions $b^2-ac\ne 0$ and also $4b^2+(a-c)^2\ne 0$, that is $b$ and $a-c$ are not both 0.
This time there is no explicit requirement that $a, c$ are nonzero,
indeed $a=c=0, b\ne 0$ gives a 2-$_V\mathcal{K}$-determined germ.

We can however reduce to two alternative normal forms, as in Table~\ref{table}, as follows.  Applying
the four vector fields~(\ref{eq:vfs}) to $dH$ the quadratic form $\phi(u,v)=au^2+2buv+cv^2$ can be changed to any
linear combination of $\phi$ and $\psi(u,v)=u\phi_v-v\phi_u=bu^2+(c-a)uv-bv^2,$ provided the conditions above are not
violated. Using $b\phi+c\psi$ we can obtain $k u^2 + uv$ for some $k$, provided $2b^2+c(c-a)\ne 0$, and using $b\phi-a\psi$ we can
obtain $uv + kv^2$ for some $k$ provided $2b^2-a(c-a)\ne 0$.  If both these reductions fail then it is easy to check that $a=c$ and $b=0$
which violates the original condition on $\phi$.

\begin{rem}
{\rm
We do not know whether this remaining smooth modulus $k$ has any geometrical
significance.  However, as in the real inflexion case (Remark~\ref{rem:mod1}), a smooth modulus is to be expected
in view of the presence of four concurrent lines in the intersection of the cone $V$ and the
zero-set of the map $(u,v,w)\mapsto (u+w, k u^2 + uv)$, to take one of the above alternatives.
Setting $u+w=z,$ the equation $uw=v^2$ becomes $u(z-u)=v^2$ and setting $z=0$ we have four
lines in this plane, $u^2+v^2=0$ and $u(ku+v)=0$.  Of course the first pair of these lines are never real.
}
\end{rem}

A $_V\mathcal{K}$ versal unfolding is given by
$(u+w, \   ku^2+uv + \lambda u + \mu v + \nu u^2)$ or $(u+w, \   uv + kv^2+\lambda u + \mu v + \nu v^2)$, where $k$ is a smooth modulus.  There
are no restrictions on the value of $k$; in particular it can be 0.
The $_V\mathcal{K}$ bifurcation set $\BB$ in this case consists of the origin only in the $(\lambda,\mu)$-plane since
$u+w=0$ is possible only for $x=y=0$,  hence $u=v=w=0$.

In the geometrical case we require $g_{31}^2-4g_{40}g_{22}\ne 0$, and $g_{31}, \ f_{20}g_{22}-f_{02}g_{40}$ are not both zero.
For $p,q$ in the family of contact maps to versally unfold the singularity we require\footnote{This is not
the same condition as that in~\cite[Prop.7.9, p.224]{IRRT} which in our notation becomes
$f_{02}(3g_{30}g_{12}-g_{21}^2)+f_{20}(3g_{21}g_{03}-g_{12}^2)\ne0$.} $g_{21}^2-3g_{12}g_{30}\ne 0$.
The inflexion points of imaginary type are isolated points of the parabolic set of $M$. They also lie on the
curve on $M$ defined by the vanishing of the {\em normal curvature} $\kappa$ of $M$. This is the same as saying
that the curvature ellipse collapses to a segment (and so has zero area). See~\cite[pp.\ 9, 17]{B-G}.
Points of the $\kappa=0$ curve on $M$ other than the inflexions of imaginary type are not distinguished by
the family of reflexion maps since in general $\kappa=0$ is not an affine invariant of $M$.

\begin{prop}
At an inflexion point of imaginary type on $M$, with 2-jet of $(f,g)$ equal to
$(f_{20}x^2 +f_{02}y^2, \ 0)$, where $f_{20}f_{02}>0$,
 the $_V\mathcal{K}$ bifurcation set consists of
the point only. A generic point of this kind is an isolated point of the parabolic set of $M$. The conditions are $g_{31}^2-4g_{40}g_{22}\ne 0$,
 $g_{22}, \ g_{40}-g_{22}$ are not both zero and $g_{21}^2-3g_{12}g_{30}\ne 0$.
\end{prop}

\section{Examples}\label{s:examples}

In this section we show how to calculate special parabolic points in practice over a whole surface $M$ given in Monge form.

A good source of examples where something interesting is happening is~\cite[pp.189-90]{B-T}.  In these examples the parabolic set undergoes
a transition as $M$ changes in a 1-parameter family,
so that a loop appears (either an elliptic island in a hyperbolic sea or vice versa), or a crossing on the parabolic set separates in a Morse
transition.  In fact from our point of view the examples of~\cite{B-T} are not quite generic since at special parabolic points, when these exist, our family
of contact maps does not versally unfold the singularity according to the criterion of \S\ref{ss:specialparabolic}.  However this is easily remedied by
additing an extra term to one of the defining equations.

For us it is not generic for a crossing or isolated point on the parabolic set to be in addition a special parabolic point, since special parabolic points
are isolated on the parabolic set. Thus when a loop of parabolic points appears on $M$ the loop will generically have no special parabolic points on it
but these can develop  as the loop expands, as the examples show.  We can check numerically that the sign of the parabolic curve, in the sense of
\S\ref{ss:parabolic}, changes at a special point, and we can calculate the type of the special point, as defined in \S\ref{ss:specialparabolic}.

\begin{figure}[!ht]
\begin{center}
\hspace*{-2cm}
\includegraphics[width=6in]{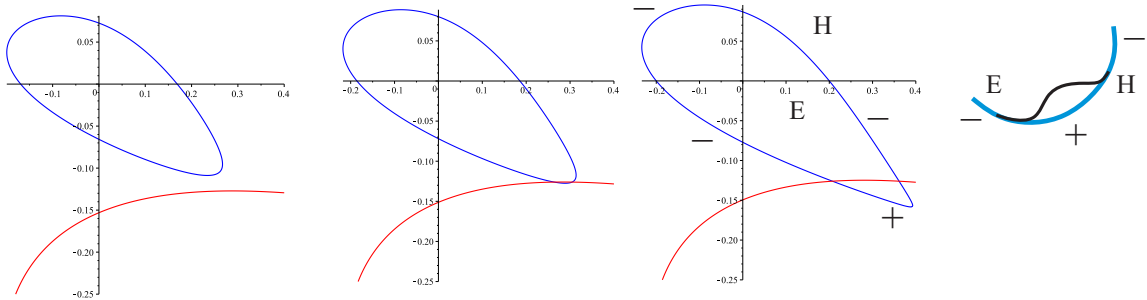}
\end{center}
\vspace*{-0.2in}
\caption{\small  The parameter plane of the curve of Example~\ref{ex1} near the origin $x=y=0$ for, left to right, $\mu = -\frac{1}{35}$,  $\mu=-\frac{1}{29}, \ \mu=-\frac{1}{25}$.  The figure shows
a loop on the parabolic set and the additional curve whose intersections with the parabolic set give special points, as in \S\ref{ss:formulas},~\S\ref{s:app}.
Two special points appear at about $\mu=-\frac{1}{29}$. The signs of the parabolic set arcs are marked in the third figure and the elliptic region E and the hyperbolic
region H. The right-hand figure is a schematic representation of the germ of a ``semi-lips'' which joins the two bifurcation sets of the special parabolic  points
immediately after their creation.  Note that this is consistent with Figure~\ref{fig:bifset-specialparabolic-realinflexion}(a) with the $-w^3$ sign. }
\label{fig:ex1}
\end{figure}

\begin{exam}\label{ex1}{\rm
Consider the family of surfaces given in Monge form by \\
$f(x,y)=xy+y^3, \ \ g(x,y)=x^2+x^2y^2+xy^3-\textstyle{\frac{1}{2}}\displaystyle y^4+\textstyle{\frac{1}{30}}\displaystyle y^5+\mu y^2,$
where the term in $y^5$ is added to the formula in~\cite[p.189]{B-T} (with $\lambda=-\frac{1}{2}$) to make the special points
generic from the family of reflexion maps, and small negative values of the parameter $\mu$ produce a loop on the parabolic set. Figure~\ref{fig:ex1} illustrates the formation
of two special points on the parabolic set as $\mu$ becomes more negative.
}
\end{exam}
We can calculate the type of the special points, and the sign of the parabolic curve on either side of them, as follows, where the calculations are
necessarily numerical rather than exact.  Having calculated numerically the parameter values $(p,q)$ of a special point, that is where the two curves
in Figure~\ref{fig:ex1} intersect we `move the origin' to this point.  This re-parametrizes $M$ near $(f(p,q), g(p,q), p, q)$ as the set of points
$ (f(x'+p, \ y'+q)-f(p,q), \ g(x'+p, \ y'+q) - g(p,q), \ x', \ y')$ where $(x', y')$ are the new coordinates in the parameter plane, with origin at $x=p, y=q$.
We can now proceed to reduce the quadratic terms of this parametrization to $(x'y', x'^2)$, ignoring any linear terms which can be removed
by a global affine transformation of $\RR^4$. Having done this, we can apply the formulas of \S\ref{ss:specialparabolic} to determine the type
of special parabolic point (elliptic or hyperbolic) and to check that it is nondegenerate and that the family of contact maps is versally unfolded.
All these calculations are straightforward and were performed in MAPLE. The same method can be used at an ordinary parabolic point to
determine whether it is positive or negative in the sense of \S\ref{ss:parabolic}.

For the example above we find that the special parabolic points are both elliptic, that is the germ of the bifurcation set is inside the elliptic
island of $M$.  We find that after reduction of the quadratic terms of $f,g$ the conditions $g_{04}=0, \ 4g_{20}g_{06}-g_{13}^2 < 0, \
5g_{12}g_{05}-3g_{03}g_{14} \ne 0$ in the notation of \S\ref{ss:specialparabolic}, all hold at both special points.  The latter condition does not
hold without the addition of the term in $y^5$ to $g$.

We also find that the sign of the parabolic points on the loop is negative for small $\mu$ before the special points appear; this is to be
expected since the sign of $y^4$ in $g(x,y)$ is $<0$.  The arc of the parabolic set between the special points consists of positive parabolic points.

\begin{figure}[!ht]
\begin{center}
\includegraphics[width=4in]{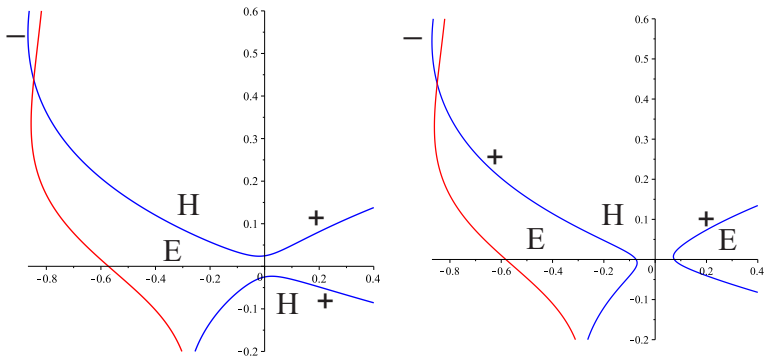}
\end{center}
\vspace*{-0.2in}
\caption{\small The parabolic set in the parameter plane for Example~\ref{ex2}, with (left) $\mu=-\frac{1}{60}$ and (right) $\mu= \frac{1}{60}$. The special
parabolic points where the two curves meet are of elliptic type; H stands for the hyperbolic region, E for the elliptic region and $+,  -$ refer to the
sign of these sections of the parabolic set, computed using the method explained above.  }
\label{fig:ex2}
\end{figure}

\begin{exam}\label{ex2}
{\rm
A second example, also adapted from~\cite{B-T}, is provided by $f(x,y)=xy+y^3, \ g(x,y)=x^2+x^2y -3x^2y^2+3y^4+y^5+\mu y^2$.
See Figure~\ref{fig:ex2} for an illustration. Calculation as above stows that the special parabolic point is elliptic and is versally unfolded by the
family of contact maps so that the bifurcation set is as described in \S\ref{ss:specialparabolic}.
Also, the signs of the parabolic set are as in the figure. Note that this transition on the parabolic set via a self-crossing is not to be confused with
the inflexion point of real type as in \S\ref{ss:realinflexion}.
}
\end{exam}

\section{Concluding remarks}\label{s:qu}

We have shown how the family of contact maps by reflexion in points of a surface $M$ in $\RR^4$ identifies the parabolic set of
$M$ and also some special but still smooth points of the parabolic set which are not part of the flat geometry of $M$ but are
affine invariants of $M$. We do not know of a different characterisation of these points.

In \cite{G-J1} it was possible to extend the bifurcation set of the family of contact maps on a surface $M$ in $\RR^3$ to
a global bifurcation set, even though it was not entirely clear what geometrical significance this had away from
the parabolic set on $M$.  In the present case, for $M$ a surface in $\RR^4$, we do not know of any reasonable
way to make the bifurcation set global.

Because of the sign attached to points of the parabolic set which changes at special parabolic points and also at self-crossings of
the parabolic set, it is possible to formulate some statements about the numbers of special points. For instance, on a smooth
closed loop of the parabolic set there must be an even number of special parabolic points (possibly zero).  Similarly on a figure-eight
component of the parabolic set there must be the same parity of special parabolic points on each loop.

It is possible in principle to extend the explicit calculations of special parabolic points, as in \S\ref{s:examples}, to the case when
the surface is parametrized in a general way, as $(A(x,y), B(x,y), C(x,y), D(x,y))$.  However there is a significant difficulty in writing down
the contact map, as in~(\ref{eq:contact-family}) which is valid for the case $C(x,y)=x, D(x,y)=y$, without an expression for $M$ as the
zero set of a submersion $\RR^4\to\RR^2$. We need to construct the contact map from parametrizations of both $M$ and its
reflexion $M^*$ in a point of $M$.  Extension to a general parametrization would allow us to examine examples such as those in~\cite{carmen}.
Even more challenging is the explicit calculation of the contact map for a surface which is given in implicit form as the zero set of
a submersion.

\appendix
\section{The additional formula for the locus of special parabolic points}\label{s:app}
Consider a surface in Monge form $(f(x,y), g(x,y),x,y)$.  For our purposes it does not matter whether $f, g$ have linear terms since they can be removed
by a global affine transformation of $\RR^4$ which will not affect the parabolic curves or special parabolic points.
The additional condition, besides (\ref{eq:parab-locus}), for a point with parameters $(p,q)$ to be a special parabolic point, is as follows.  We use the notation
of \S\ref{ss:formulas}.

Let
\[ \Theta_1=f_{12}g_{22}-f_{22}g_{12}, \  \Theta_2=f_{11}g_{22}-f_{22}g_{11}, \ \Theta_3 = f_{11}g_{12}-f_{12}g_{11} \]
\[ \Phi_1=f_{11}g_{11}g_{22}-2f_{11}g_{12}^2+2f_{12}g_{11}g_{12}-f_{22}g_{11}^2\]
\[ \Phi_2=f_{11}g_{11}f_{22}-2f_{12}^2g_{11}+2f_{11}f_{12}g_{12}-f_{11}^2g_{22}\]

Then the condition is
\[ \Theta_1^2\Phi_1f_{1111} -2\Theta_1\Theta_2\Phi_1 f_{1112} +6\Theta_1\Theta_3\Phi_1 f_{1122} -2\Theta_2\Theta_3\Phi_1 f_{1222} \]

\vspace*{-0.4in}

\[ +\Theta_3^2\Phi_1 f_{2222}+ \Theta_1^2\Phi_2g_{1111}-2\Theta_1\Theta_2\Phi_2g_{1112} +6\Theta_1\Theta_3\Phi_2g_{1122}\]

\vspace*{-0.25in}

\[-2\Theta_2\Theta_3\Phi_2g_{1222} + \Theta_3^2\Phi_2g_{2222}=0.  \]

In the case that $f_{11}=0, f_{12}\ne 0, f_{22}=0, g_{11}\ne 0, g_{12}=0, g_{22}=0$ this reduces to $g_{2222} = 0$, as we expect from \S\ref{ss:specialparabolic} where
the condition appears as $g_{04}=0$ when we are working at the origin.

\bigskip\noindent
{\sc Acknowledgements}  We are very grateful to Victor Goryunov for helpful conversations, in particular for suggesting the
correct definition of the bifurcation set in our case. Also, the first two authors acknowledge
the hospitality and support of the University of S\~{a}o Paulo at S\~{a}o Carlos in 2014; the first author acknowledges the
hospitality and support of the Center for Advanced Studies at the Technical University of Warsaw; and the second and third
authors acknowledge the hospitality and support of the University of Liverpool, the London Mathematical Society and the UK research council EPSRC during a visit in 2016.

\noindent
Peter J Giblin, Department of Mathematical Sciences, The University of Liverpool, Liverpool L69~7ZL, UK, pjgiblin@liv.ac.uk

\smallskip\noindent
Stanis\l aw Janeczko, Faculty of Mathematics and Information Science, Warsaw University of Technology, Pl. Politechniki 1, 00-661 Warsaw, Poland; and
Institute of Mathematics Polish Academy of Sciences, Sniadeckich 
8, 00-656 Warsaw, Poland.
 S.Janeczko@mini.pw.edu.pl

\smallskip\noindent
Maria Aparecida Ruas, ICMC, USP, Avenida Trabalhador S\~{a}o-Carlense, 400  Centro, CEP: 13566-590, S\~{a}o Carlos SP, Brazil,
maasruas@icmc.usp.br

\end{document}